\documentclass[11pt]{article}
\usepackage{amsfonts,amsmath,amsxtra}
\usepackage{latexsym}
\usepackage{amssymb}

\def\hybrid{\topmargin 0pt      \oddsidemargin 0pt
        \headheight 0pt \headsep 0pt
        \textwidth 16.5cm
        \textheight 23cm
        \marginparwidth 0.0in
        \parskip 5pt plus 1pt   \jot = 1.5ex}
\catcode`\@=11
\def\marginnote#1{}
\newcount\hour
\newcount\minute
\newtoks\amorpm
\hour=\time\divide\hour by60 \minute=\time{\multiply\hour by60
\global\advance\minute by-\hour}
\edef\standardtime{{\ifnum\hour<12 \global\amorpm={am}%
        \else\global\amorpm={pm}\advance\hour by-12 \fi
        \ifnum\hour=0 \hour=12 \fi
      \number\hour:\ifnum\minute<10 0\fi\number\minute\the\amorpm}}
\edef\militarytime{\number\hour:\ifnum\minute<10 0\fi\number\minute}

\def\draftlabel#1{{\@bsphack\if@filesw {\let\thepage\relax
   \xdef\@gtempa{\write\@auxout{\string
      \newlabel{#1}{{\@currentlabel}{\thepage}}}}}\@gtempa
   \if@nobreak \ifvmode\nobreak\fi\fi\fi\@esphack}
        \gdef\@eqnlabel{#1}}
\def\@eqnlabel{}
\def\@vacuum{}
\def\draftmarginnote#1{\marginpar{\raggedright\scriptsize\tt#1}}

\def\draft{\oddsidemargin -0.1truein
        \def\@oddfoot{\sl preliminary draft \hfil
        \rm\thepage\hfil\sl\today\quad\militarytime}
        \let\@evenfoot\@oddfoot \overfullrule 3pt
        \let\label=\draftlabel
        \let\marginnote=\draftmarginnote
\def\@eqnnum{{\rm (\theequation)}
\rlap{\kern\marginparsep\tt\@eqnlabel}%
\global\let\@eqnlabel\@vacuum}  }

\newfont{\Bbbb}{msbm7 scaled 1\@ptsize00}
\newcommand{\zs}{\raise-1pt\hbox{$\mbox{\Bbbb Z}$}}

scaled 1\@ptsize00

\font\sevenmsa=msam6 
\newfam\msafam
\textfont\msafam=\sevenmsa
\def\hexnumber@#1{\ifnum#1<10 \number#1\else
\ifnum#1=10 A\else\ifnum#1=11 B\else\ifnum#1=12 C\else \ifnum#1=13
D\else\ifnum#1=14 E\else\ifnum#1=15 F\fi\fi\fi\fi\fi\fi\fi}
\def\msa@{\hexnumber@\msafam}
\def\llcorner{\delimiter"4\msa@78\msa@78 }
\def\lrcorner{\delimiter"5\msa@79\msa@79 }
\mathchardef\blacktriangleright="3\msa@49
\mathchardef\blacktriangleleft="3\msa@4A \font\tenmsb=msbm10 scaled
1\@ptsize00
\newfam\msbfam
\textfont\msbfam=\tenmsb \scriptfont\msbfam=\tenmsb

\newdimen\Squaresize \Squaresize=14pt
\newdimen\Thickness \Thickness=0.5pt

\def\Square#1{\hbox{\vrule width \Thickness
   \vbox to \Squaresize{\hrule height \Thickness\vss
      \hbox to \Squaresize{\hss#1\hss}
   \vss\hrule height\Thickness}
\unskip\vrule width \Thickness} \kern-\Thickness}

\def\Vsquare#1{\vbox{\Square{$#1$}}\kern-\Thickness}

\def\numberbysection{\@addtoreset{equation}{section}
        \def\theequation{\thesection.\arabic{equation}}}
\numberbysection

\renewcommand{\theequation}{\thesection.\arabic{equation}}
\def\titlepage{\@restonecolfalse\if@twocolumn\@restonecoltrue\onecolumn
     \else \newpage \fi \thispagestyle{empty}\c@page\z@
        \def\thefootnote{\fnsymbol{footnote}} }

\def\endtitlepage{\if@restonecol\twocolumn \else  \fi
        \def\thefootnote{\arabic{footnote}}
        \setcounter{footnote}{0}}  
\relax

\hybrid
\makeatletter
\newdimen\normalarrayskip            
\newdimen\minarrayskip               
\normalarrayskip\baselineskip \minarrayskip\jot
\newif\ifold             \oldtrue            \def\new{\oldfalse}
\def\arraymode{\ifold\relax\else\displaystyle\fi}
\def\eqnumphantom{\phantom{(\theequation)}} 
\def\@arrayskip{\ifold\baselineskip\z@\lineskip\z@
     \else
     \baselineskip\minarrayskip\lineskip1\baselineskip\fi}

\def\@arrayclassz{\ifcase \@lastchclass \@acolampacol \or
\@ampacol \or \or \or \@addamp \or
   \@acolampacol \or \@firstampfalse \@acol \fi
\edef\@preamble{\@preamble
  \ifcase \@chnum
     \hfil$\relax\arraymode\@sharp$\hfil
     \or $\relax\arraymode\@sharp$\hfil
     \or \hfil$\relax\arraymode\@sharp$\fi}}

\def\@array[#1]#2{\setbox\@arstrutbox=\hbox{\vrule
     height\arraystretch \ht\strutbox
     depth\arraystretch \dp\strutbox
width\z@}\@mkpream{#2}\edef\@preamble{\halign \noexpand\@halignto
\bgroup \tabskip\z@ \@arstrut \@preamble \tabskip\z@ \cr}%
\let\@startpbox\@@startpbox \let\@endpbox\@@endpbox
  \if #1t\vtop \else \if#1b\vbox \else \vcenter \fi\fi
  \bgroup \let\par\relax
  \let\@sharp##\let\protect\relax
  \@arrayskip\@preamble}
\def\eqnarray{\stepcounter{equation}%
              \let\@currentlabel=\theequation
              \global\@eqnswtrue
              \global\@eqcnt\z@
              \tabskip\@centering              
              \let\\=\@eqncr
              $$%
            \halign to \displaywidth  \bgroup
             \eqnumphantom \@eqnsel
      \hskip\@centering                               
    $\displaystyle  \tabskip\z@ {##}$%
    &\global\@eqcnt\@ne \hskip 2\arraycolsep
         $ \displaystyle  \arraymode{##}$\hfil
    &\global\@eqcnt\tw@ \hskip 2\arraycolsep
         $\displaystyle\tabskip\z@{##}$\hfil
         \tabskip\@centering
    &{##}\tabskip\z@\cr}
\makeatother
\newcommand{\RR}{{\mathbb{R}}}

\newcommand{\ZZ}{{\mathbb{Z}}}

\newcommand{\C}{\mathbb{C}}

\def\IC{\mathbb{C}}

\def\IZ{\mathbb{Z}}

\def\CH {\mathcal{H}}

\def\CJ {\mathcal{J}}

\def\CP {\mathcal{P}}
\def\CQ {\mathcal{Q}}

\def\ch{{\cal H}}

\def\l {{\lambda}}

\def\l{\lambda}

\def\pr {\partial}

\def\Tr{{\rm Tr}}

\def\ch{{\rm ch}}

\def\frak{\mathfrak}

\newtheorem{te}{Theorem}[section]
\newtheorem{de}{Definition}[section]
\newtheorem{prop}{Proposition}[section]           
\newtheorem{cor}{Corollary}[section]
\newtheorem{lem}{Lemma}[section]
\newtheorem{ex}{Example}[section]
\newtheorem{rem}{Remark}[section]
\newcommand\bqa{\begin{eqnarray}}
\newcommand\eqa{\end{eqnarray}}
\def\be{\begin{eqnarray}\new\begin{array}{cc}}
\def\ee{\end{array}\end{eqnarray}}

\def\beq{\begin{equation}}
\def\eeq{\end{equation}}
\def\bse{\begin{subequations}}                
\def\ese{\end{subequations}}
\def\bp{\begin{pmatrix}}
\def\ep{\end{pmatrix}}

\def\i{\imath}

\def\stack#1#2{\raise0.7pt\hbox{$\mathrel{\mathop{#2}\limits^{#1}}$}}
\def\tr{\triangleright}
\def\tl{\triangleleft}
\def\sem{\mathsurround=0pt \raise1pt
\hbox{$\scriptscriptstyle>\!\!$}\:\!\!\tl}
\def\mes{\mathsurround=0pt \tr\!\:\!\raise0.8pt
\hbox{$\scriptscriptstyle\!\!<$}\,}
\def\]{\mathsurround=0pt ]\raise-2pt\hbox{$_\ast$}}

\def\<{\langle}
\def\>{\rangle}

\def\CQ{{\cal Q}}
\def\frak{\mathfrak}

\def\ch{{\cal H}}

\def\CH{\mathcal{H}}

\def\we{\raise-1pt\hbox{$\,\stackrel{\wedge}{,}\,$}}
\def\tr{{\rm tr}\,}
\def\Tr{{\rm Tr}\,}
\def\pr {\partial}

\newcounter{pac}[section]

0

0

\title{\bf On $q$-deformed $\mathfrak{gl}_{\ell+1}$-Whittaker function I}
\begin{document}
\author{Anton Gerasimov, Dimitri Lebedev, and Sergey Oblezin}
\date{}

\maketitle


\renewcommand{\abstractname}{}

\begin{abstract}
\noindent {\bf Abstract}. We propose  new explicit form of
$q$-deformed Whittaker functions solving
$q$-deformed $\mathfrak{gl}_{\ell+1}$-Toda chains.
In the limit $q\to 1$ constructed solutions reduce
to classical class one $\mathfrak{gl}_{\ell+1}$-Whittaker
functions in the form proposed by Givental.
 An important  property of the proposed
expression for the $q$-deformed $\mathfrak{gl}_{\ell+1}$-Whittaker
function is that it can be represented as a character of $\IC^*\times
GL(\ell+1)$. This provides a $q$-version of the
Shintani-Casselman-Shalika formula for  $p$-adic Whittaker function.
 The Shintani-Casselman-Shalika  formula is recovered in the limit
 $q\to 0$ when the $q$-deformed Whittaker function
is reduced to a character of a finite-dimensional representation
of $\mathfrak{gl}_{\ell+1}$ expressed through Gelfand-Zetlin
bases.

\end{abstract}

\section*{Introduction}

Whittaker functions corresponding
to semisimple finite-dimensional Lie algebras   arise in
various parts of modern mathematics. In particular,  these functions appear
in representation theory as matrix
elements of infinite-dimensional representations, in the theory of
quantum integrable systems as  common eigenfunction of Toda chain
 quantum Hamiltonians, in string theory
as generating functions of correlators in Type A topological 
string theory on flag manifolds and in number theory 
in a description of local Archimedean
$L$-factors corresponding to automorphic representations.
Although much studied, Whittaker functions seems have
some deep properties that are not yet fully revealed.

In this paper we study  the $q$-deformed
$\mathfrak{gl}_{\ell+1}$-Whittaker functions.
The $q$-deformed  Whittaker function can be identified with a common eigenfunction
of a set of commuting $q$-deformed Toda chain Hamiltonians.  This
$q$-deformed Toda chain (also known as the relativistic Toda chain \cite{Ru}) was
discussed in terms of representation theory of  quantum groups  in
\cite{Se1}, \cite{Et}, \cite{Se2} and  an integral representation for
the $q$-deformed $\mathfrak{gl}_{\ell+1}$-Whittaker function was
constructed in \cite{KLS}. Recently the $q$-deformed Toda chain
attracts special interest due to its connection with quantum
$K$-theory of flag manifolds \cite{GiL}.
 In this paper we pursue another direction.  Our principal motivation
 to study $q$-{\it deformed}  Whittaker functions is that in this, more general
setting, some important hidden properties of classical Whittaker
functions become visible. 

 The main result of the paper is given by Theorem \ref{mainTH} where a
new expression for the $q$-deformed $\mathfrak{gl}_{\ell+1}$-Whittaker
function ( for $q<1$) is introduced. As a simple corollary of 
Theorem \ref{mainTH}, the $q$-deformed $\mathfrak{gl}_{\ell+1}$-Whittaker
function can be represented as a character of $\IC^*\times
GL(\ell+1)$. In the limit $q\to 1$ this leads to a similar
representation of classical $\mathfrak{gl}_{\ell+1}$-Whittaker
function. This representation is not easy to perceive looking
directly at the classical Whittaker functions.
The importance of this representation of ($q$-deformed)
$\mathfrak{gl}_{\ell+1}$-Whittaker function becomes obvious if we
notice that in the limit $q\to 0$ the constructed $q$-deformed
Whittaker function reduces to $p$-adic Whittaker function. In this
limit the representation as a character reduces to
well-known Shintani-Casselman-Shalika representation of $p$-adic
$GL_{\ell+1}$-Whittaker function as a character of a
finite-dimensional representation of $GL_{\ell+1}$ 
\cite{Sh},\cite{CS}. Thus the representation of ($q$-deformed)
$\mathfrak{gl}_{\ell+1}$-Whittaker function as a character can
be considered as a $q$-version  of Shintani-Casselman-Shalika
representation. Indeed, the constructed $q$-deformed
Whittaker function is equal to zero outside a dominant
weight cone of $\mathfrak{gl}_{\ell+1}$ similarly to the
Shintani-Casselman-Shalika $p$-adic Whittaker function.

We expect that the representation of the classical Whittaker
function  as a character should provide important  insights
into the arithmetic  geometry at an infinite  place of $\overline{{\rm
Spec}(\IZ)}$. Let us also remark that taking into account results
\cite{CS} one should expect that in the case of an arbitrary
semisimple Lie algebra $\mathfrak{g}$, $q$-deformed
$\mathfrak{g}$-Whittaker function should be  given by  a character
of $\IC^*\times {}^LG(\IC)$  where ${\rm Lie}({}^LG)={}^L\mathfrak{g}$ is a
Langlands dual Lie algebra.

It is worth mentioning that the $q\to 1$  limit of the explicit
expression of the $q$-deformed Whittaker function proposed
in this paper reduces to the integral representations for classical
Whittaker functions introduced by Givental 
\cite{Gi},\cite{GKLO}. We consider this as
a sign of an ``arithmetic nature'' of this integral representation.
On the other hand the explicit  solution has an
obvious relation with Gelfand-Zetlin parametrization of
finite-dimensional representations of $\mathfrak{gl}_{\ell+1}$ (and
precisely reproduces Gelfand-Zetlin form of characters of
finite-dimensional representations in the limit $q\to 0$). This duality
of Gelfand-Zetlin and Givental representations was already noticed
in \cite{GLO}.

Let us comment on our approach to derivation of  explicit
expressions for $q$-deformed Whittaker functions. It is known
\cite{Et} that defining difference equations for Macdonald
polynomials are transformed into  $q$-deformed Toda chain
eigenfunction equations in a certain limit. This is a simple
generalization of the Inozemtsev limit \cite{I} transforming
Calogero-Sutherland integrable model into standard Toda chain. The
other ingredient we use is a recursive construction of Macdonald
polynomials (analogous to  the recursive construction  for
($q$-deformed) Toda chain  eigenfunctions \cite{KL1}, \cite{KLS}).
We combine these results to obtain recursive expression for
$q$-deformed $\mathfrak{gl}_{\ell+1}$-Whittaker functions satisfying
$q$-deformed $\mathfrak{gl}_{\ell+1}$-Toda chain eigenfunction
equations.

The explicit form of the $q$-deformed Whittaker function implies
various interesting interpretations. This 
includes connections with representation theory
(via characters of Demazure modules), geometry of quiver varieties, 
quantum cohomology of flag manifolds  and 
will be discussed elsewhere \cite{GLO2}.

Finally note that eigenfunctions of
$q$-deformed Toda chain   were  discussed previously (e.g.
\cite{KLS},\cite{GKL1},\cite{BF} and \cite{FFJMM}).
The relation of these constructions
with the one proposed in this paper is an interesting
question which deserves further considerations.

The paper is organized as follows. In Section 1 we recall  a systems
of mutually commuting difference Macdonald-Ruijsenaars operators and
recursive construction of their common eigenfunctions. In Section 2
we derive recursive expression for  solutions of $q$-deformed
$\mathfrak{gl}_{\ell+1}$-Toda chain. In Section 3 various limiting
cases elucidating the construction of  the $q$-deformed
$\mathfrak{gl}_{\ell+1}$-Whittaker functions are discussed. In
Section 4 details of the proof of the Theorem \ref{mainTH} are
given.

{\em Acknowledgments}: The research of AG was  partly supported by
SFI Research Frontier Programme and Marie Curie RTN Forces Universe
from EU. The research of SO is partially supported by  RF President
Grant MK-134.2007.1.

\section{Macdonald-Ruijsenaars difference operators }

In this section we recall  relevant facts
from the theory of Macdonald polynomials  ( see e.g.
 \cite{Mac}, \cite{Kir}, \cite{AOS}).

Consider symmetric polynomials in  variables $(x_1,\ldots,x_{\ell+1})$ over the
field $\mathbb{Q}(q,t)$ of rational functions in $q,t$.
Given a partition $\Lambda=(0\leq\Lambda_1\leq\Lambda_2\leq\ldots\leq\Lambda_{\ell+1})$,
denote by the same symbol $\Lambda$ the Young diagram containing $\ell+1$
rows with $\Lambda_k$ boxes in the $k$-th row; and the upper row
having  the maximal length $\Lambda_{\ell+1}$. 

Let $m_\Lambda$ and $\pi_\Lambda$ be  polynomial basises of the
space of symmetric polynomials indexed by partitions $\Lambda$:
$$
m_\Lambda=\sum_{\sigma\in\mathfrak{S}_{\ell+1}}\,
x_{\sigma(1)}^{\Lambda_1}x_{\sigma(2)}^{\Lambda_2}\cdot \ldots\cdot
x_{\sigma(\ell+1)}^{\Lambda_\ell+1},
$$
$$
\pi_\Lambda\,=\,\pi_{\Lambda_1}\pi_{\Lambda_2}\cdot\ldots\cdot
\pi_{\Lambda_{\ell+1}},\hspace{1.5cm}\pi_n=\sum_{k=1}^{\ell+1}\,x_k^{n},
$$
where $\mathfrak{S}_{\ell+1}$ is the permutation group. Define a
scalar product $\<\,,\,\>_{q,t}$ on the space of symmetric functions
over $\mathbb{Q}(q,t)$ as follows
$$
\<\pi_\Lambda,\,\pi_{\Lambda'}\>_{q,t}\,=\,\delta_{\Lambda,\Lambda'}\cdot
z_\Lambda(q,t),
$$
where
$$
z_\Lambda(q,t)\,=\,\prod_{n\geq1}\,n^{m_n}m_n!\cdot \prod_{k=1}^N
\frac{1-q^{\Lambda_k}}{1-t^{\Lambda_k}},\hspace{1.5cm}
m_n=\bigl|\{k|\,\Lambda_k=n\}\bigr|.
$$
In the following we always imply $q<1$.

\begin{de} Macdonald polynomials
$P^{\mathfrak{gl}_{\ell+1}}_\Lambda=
P^{\mathfrak{gl}_{\ell+1}}_{\Lambda}(x;q,t)$
are symmetric polynomial  function over $\mathbb{Q}(q,t)$ such that\\
$$
P^{\mathfrak{gl}_{\ell+1}}_\Lambda\,=\,
m_\Lambda+\sum_{\Lambda'<\Lambda} u_{\Lambda\Lambda'}m_{\Lambda'},
$$
with $u_{\Lambda\Lambda'}\in\mathbb{Q}(q,t)$, and
 for $\Lambda\neq\Lambda'$
$$
\bigl\<P^{\mathfrak{gl}_{\ell+1}}_\Lambda,\,
P^{\mathfrak{gl}_{\ell+1}}_{\Lambda'}\bigr\>_{q,t}\,=\,0.
$$
\end{de}

Macdonald polynomials are eigenfunctions of a set of mutually
commuting  Macdonald-\\Ruijsenaars difference operators \cite{Mac},
\cite{Ru} \be\label{macdif}
H^{\mathfrak{gl}_{\ell+1}}_{r}=\sum_{I_r} \, t^{r(r-1)/2}\,
\prod_{i\in I_r,\,j\notin I_r}
\frac{tx_i-x_j}{x_i-x_j}\,\,\,\,\prod_{m\in I_r}T_m,  \qquad
r=1,\ldots,\ell+1,\ee where  the sum is over ordered subsets
$$
I_r=\{i_1<i_2<\ldots<i_r\}\subset\{1,2,\ldots,\,\ell+1\}.
$$
The simplest operator of this kind is given by
$$
H_1^{\mathfrak{gl}_{\ell+1}}=\sum_{i=1}^{\ell+1}\,\prod_{j,\,j\neq
i} \frac{tx_i-x_j}{x_i-x_j}\,q^{x_i\pr_{x_i}}.
$$
The eigenvalues of $H_r$ are given by (see e.g. \cite{EK})
$$
H^{\mathfrak{gl}_{\ell+1}}_{r}\,
P^{\mathfrak{gl}_{\ell+1}}_{\Lambda}(x;q,t)=c_{\Lambda}^{r}\,
P^{\mathfrak{gl}_{\ell+1}}_{\Lambda}(x;q,t),
$$
$$
c_{\Lambda}^r=
\chi_r\bigl(q^{\sum_{i=1}^{\ell+1}E_{ii}\Lambda_i}t^{\sum_{i=1}^{\ell+1}
E_{i,i}(\ell+1-i)}\bigr)\,=\,\sum_{I_r}
\prod_{i\in I_r}\,q^{\Lambda_{i}}\,\, t^{\ell+1-i},
$$
where $E_{i,j}$ are standard generators of $\mathfrak{gl}_{\ell+1}$,
 $I_r=(i_1<i_2<\ldots<i_r)\subset\{1,2,\ldots,\ell+1\}$
and $\chi_r(g)$ are the character of fundamental representations
$V_r=\bigwedge^{r}\C^{\ell+1}$ of $\mathfrak{gl}_{\ell+1}$.
In terms of a generating series
$$
H^{\mathfrak{gl}_{\ell+1}}(\xi)=\sum_{r=0}^{\ell+1} \xi^{\ell+1-r}
H^{\mathfrak{gl}_{\ell+1}}_r,\qquad H_0=1
$$
we have
$$
 H^{\mathfrak{gl}_{\ell+1}}(\xi)\,
P^{\mathfrak{gl}_{\ell+1}}_{\Lambda}(x;q,t)=
\prod_{i=1}^{\ell+1}(\xi+t^{\ell+1-i} q^{\Lambda_i})
\,P^{\mathfrak{gl}_{\ell+1}}_{\Lambda}(x;q,t).
$$
Define  the following scalar product on symmetric functions of
$(\ell+1)$-variables $x_1,\ldots ,x_{\ell+1}$
$$
\<f,g\>'_{q,t}=\frac{1}{(\ell+1)!}\oint_{x_1=0}\cdots \oint_{x_{\ell+1}=0}
\,\prod_{i=1}^{\ell+1}\,\frac{dx_i}{2\pi\imath x_i} \,
f(x^{-1})\,g(x)\,\Delta(x|q,t),
$$
where
$$
\Delta(x|q,t)=\prod_{i\neq j}\,\prod_{n=0}^{\infty}
\frac{1-x_ix_j^{-1}q^n}{1-tx_ix_j^{-1}q^n}.
$$
Then difference operators $H^{\mathfrak{gl}_{\ell+1}}_r$ are
self-adjoint with respect to $\<,\,\>'_{q,t}$:
$$
\<f,H^{\mathfrak{gl}_{\ell+1}}_{r}\,g\>'_{q,t}\,=\,
\<H^{\mathfrak{gl}_{\ell+1}}_{r}\,f,g\>'_{q,t}\,\,.
$$
For Macdonald polynomials one has an analog of Cauchy-Littlewood
formula
$$
C_{\ell+1,m+1}(x,y|q,t) =\sum_\Lambda
P_{\Lambda}(x;q,t)\,P_{\Lambda}(y;q,t)\,\,b_{\Lambda}(q,t), \qquad
m\leq\ell,
$$
where the sum is over all Young diagrams of $\mathfrak{gl}_{m+1}$ and
\be\label{mackern}
C_{\ell+1,m+1}(x,y|q,t)=\prod_{i=1}^{\ell+1}\prod_{j=1}^{m+1}
\prod_{n=0}^{\infty} \,\frac{1-tx_iy_jq^n}{1-x_iy_jq^n},
\ee
$$
b_{\Lambda}(q,t)\,=\,
\frac{1}{\<P^{\mathfrak{gl}_{\ell+1}}_\Lambda,\,
P^{\mathfrak{gl}_{\ell+1}}_\Lambda\>_{q,t}}\,=\,
\prod_{n=1}^N\prod_{k=0}^{n-1}\,\,
\prod_{m=\Lambda_{n,n-k}}^{\Lambda_{n,n-k-1}-1}\,
\frac{1-q^{m+1}t^k}{1-q^mt^{k+1}},\qquad
\Lambda_{n,j}=\Lambda_n-\Lambda_j,\,\,\,j<n.
$$

\begin{prop}\cite{AOS}
 The following relations hold

1.\be\label{macrec} P^{\mathfrak{gl}_{\ell+1}}_{\Lambda}(x;q,t)\,=\,
\frac{ \<P^{\mathfrak{gl}_\ell}_{\Lambda},
P^{\mathfrak{gl}_\ell}_{\Lambda}\>_{q,t}}{\ell!
\<P^{\mathfrak{gl}_\ell}_{\Lambda},
P^{\mathfrak{gl}_\ell}_{\Lambda}\>'_{q,t}}\,\, \\
\times \,\,\oint_{x_1=0}\cdots \oint_{x_{\ell+1}=0}
 \prod_{i=1}^{\ell}\,\frac{dy_i}{2\pi\i y_i}\,\,
C_{\ell+1,\ell}(x,y^{-1}|q,t) P^{\mathfrak{gl}_{\ell}}_{\l}(y;q,t)
\,\Delta(y|q,t). \ee

2. \be\label{macbax}  P^{\mathfrak{gl}_{\ell+1}}_{\Lambda}(x;q,t)=
\frac{\<P^{\mathfrak{gl}_{\ell+1}}_{\Lambda},
P^{\mathfrak{gl}_{\ell+1}}_{\Lambda}\>_{q,t}}{(\ell+1)!
\<P^{\mathfrak{gl}_{\ell+1}}_{\Lambda},
P^{\mathfrak{gl}_{\ell+1}}_{\Lambda}\>'_{q,t}}\,\, \\
\times \,\,\oint_{x_1=0}\cdots \oint_{x_{\ell+1}=0}
 \prod_{i=1}^{\ell+1}\,\frac{dy_i}{2\pi\i y_i}\,\,
C_{\ell+1,\ell+1}(x,y^{-1}|q,t)
P^{\mathfrak{gl}_{\ell+1}}_{\Lambda}(y;q,t) \,\Delta(y|q,t). \ee

3.\be P^{\mathfrak{gl}_{\ell+1}}_{\Lambda+(\ell+1)^k}(x;q,t)
=\Big(\prod_{j=1}^{\ell+1}\,x_j^{k}\Big)\,
P^{\mathfrak{gl}_{\ell+1}}_{\Lambda}(x;q,t). \ee Here
$\Lambda+(\ell+1)^k$ is a Young diagram obtained from $\Lambda$ by a
substitution $\Lambda_j\to \Lambda_j+k$ and
$$
\<P^{\mathfrak{gl}_{\ell+1}}_\Lambda,P^{\mathfrak{gl}_{\ell+1}}_\Lambda\>'_{q,t}\,=\,
\prod_{1\leq i<j\leq {\ell+1}}\prod_{n=0}^\infty
\frac{1-t^{j-i}q^{\Lambda_i-\Lambda_j+n}}
{1-t^{j-i+1}q^{\Lambda_i-\Lambda_j+n}}\cdot
\frac{1-t^{j-i}q^{\Lambda_i-\Lambda_j+n+1}}
{1-t^{j-i-1}q^{\Lambda_i-\Lambda_j+1}},
$$
$$
\<P^{\mathfrak{gl}_{\ell+1}}_\Lambda,\,P^{\mathfrak{gl}_{\ell+1}}_\Lambda\>_{q,t}\,=\,
\prod_{n=1}^N\prod_{k=0}^{n-1}\,\,
\prod_{m=\Lambda_{n,n-k}}^{\Lambda_{n,n-k-1}-1}\,
\frac{1-q^mt^{k+1}}{1-q^{m+1}t^k},\qquad
\Lambda_{n,j}=\Lambda_n-\Lambda_j,\,\,\,j<n.
$$
\end{prop}
These relations provide a recursive construction
of Macdonald polynomials corresponding to  arbitrary Young diagrams.

\begin{rem} Relations \eqref{macrec} and \eqref{macbax} are analogous
to the action of recursion and Baxter operators on Whittaker functions considered
in \cite{GLO}.
\end{rem}
It  follows from above considerations that the following
intertwining relations hold.  

\begin{prop}\label{itermac} Let ${H}_k^{\frak{gl}_{\ell+1}}(x)$ and
$C_{\ell+1,\ell}(x,y|q,t)$ be given by (\ref{macdif}) and
(\ref{mackern}) respectively.
 The following intertwining  relations  hold \be\label{macbraid}
{{H}}_k^{\frak{gl}_{\ell+1}}(x)C_{\ell+1,\ell}(x,y|q,t)=
(t^{k-1}{{H}}^{\frak{gl}_{\ell}}_{k-1}(y)+
t^k{{H}}^{\frak{gl}_\ell}_k(y))C_{\ell+1,\ell}(x,y|q,t), \qquad
k=1,\ldots,\ell+1. \ee
\end{prop}

\section{ $q$-deformed  $\mathfrak{gl}_{\ell+1}$-Whittaker
function}

Quantum Hamiltonians of  $q$-deformed $\mathfrak{gl}_{\ell+1}$-Toda
chain can be considered as particular degeneration of Macdonald-Ruijsenaars
difference operators \cite{Et}.
 This is an analog of the  Inozemtsev limit \cite{I}
producing Toda chains from Calogero-Sutherland models. In this
section we give an explicit expression for $q$-deformed
$\mathfrak{gl}_{\ell+1}$-Whittaker function obtained using a
degeneration of recursive operators for Macdonald polynomials
discussed in the previous section. We also consider interesting
features of the obtained expressions. Details of the proof will be
given in the last section.

Quantum $\mathfrak{gl}_{\ell+1}$-Toda chain  is defined by a  set of
$\ell+1$  mutually commuting functionally independent quantum
Hamiltonians $\CH_r^{\mathfrak{gl}_{\ell+1}}$, $r=1,\ldots ,\ell+1$
\be\label{comm} \CH_r^{\mathfrak{gl}_{\ell+1}}(x)\,=\,\sum_{I_r}\,
\bigl(X_{i_1}^{1-\delta_{i_1,\,1}}\cdot
X_{i_2}^{1-\delta_{i_2-i_1,\,1}}\cdot\ldots\cdot
X_{i_r}^{1-\delta_{i_r-i_{r-1},\,1}}\bigr) T_{i_1}\cdot\ldots\cdot
T_{i_r}, \ee where  summation goes over ordered subsets $I_r
=\{i_1<i_2<\cdots<i_r\} $ of $\{1,2,\cdots,\ell+1\}$ and
$X_i(x)\,:=\,1-x_ix_{i-1}^{-1}$, $1<i\leq\ell+1$  with
$X_1=1$, $T_i=q^{x_i\pr_{x_i}}$.  The simplest operator is given by
$$
\CH_1^{\mathfrak{gl}_{\ell+1}}(x)=T_1+\sum_{i=1}^{\ell}
(1-x_{i+1}x_i^{-1})T_{i+1}.
$$
Common eigenfunctions of the Hamiltonians are given by
$q$-deformed $\mathfrak{gl}_{\ell+1}$-Whittaker
functions  \cite{Et}
 \be\label{qToda}
\CH_r^{\mathfrak{gl}_{\ell+1}}(x)\,\,\,
\Psi^{\mathfrak{gl}_{\ell+1}}_{\l_1,\cdots ,\l_{\ell+1}}
(x_1,\ldots,x_{{\ell}+1})=(\sum_{ I_r} \prod\limits_{i\in I_r}
q^{\l_i} )\,\, \Psi^{\mathfrak{gl}_{\ell+1}}_{\l_1,\cdots
,\l_{\ell+1}} (x_1,\ldots,x_{\ell+1}),
\ee
where $\l_i$, $i=1,\ldots \ell=1$ are real numbers.

Note that the set of equations \eqref{qToda} allows an
infinite-dimensional linear space of solutions (due to a possibility
to multiply any solution on an arbitrary function $f(x_1,\ldots,
x_{\ell+1})$ periodic with respect to the shifts $x_i\to
x_i+m_i$, $m_i\in\IZ$). To obtain a finite-dimensional space of
solutions we specify the variables to the lattice
$\ZZ^{\ell+1}$ as follows
$$
x_j=q^{p_{\ell+1,j}+j-1},\hspace{1.5cm} p_{\ell+1,j}\in \IZ,
j=1,\ldots,\ell+1.
$$
We shall use the following notation
$\underline{p}_{\ell+1}=(p_{\ell+1,1},\ldots,p_{\ell+1,\ell+1})$.
The complete set of commuting Hamiltonians (\ref{comm}) can be restricted to the lattice
$\IZ^{\ell+1}$ using the substitution
$X_i(\underline{p}_{\ell+1})=1-q^{p_{\ell+1,i}-p_{\ell+1,i-1}+1}$,
$X_1(\underline{p}_{\ell+1})=1$ and
$T_if(\underline{p}_{\ell+1})=f(\underline{\widetilde{p}}_{\ell+1})$
with $\widetilde{p}_{\ell+1,k}=p_{\ell+1,k}+\delta_{k,i}$. 
Thus the first non-trivial Hamiltonian is given by:
$$
 {\ch}_1^{\mathfrak{gl}_{\ell+1}}(\underline{p}_{\ell+1})=
 T_1+\sum\limits_{i=1}^{\ell}(1-q^{p_{\ell+1,i+1}-p_{\ell+1,i}+1})T_{i+1}.
$$
 We shall be interested in the solution of the eigenvalue problem
of  $q$-deformed Toda chain on the lattice $\ZZ^{\ell+1}$:
\be\label{eiglat}
\ch_r^{\mathfrak{gl}_{\ell+1}}(\underline{p}_{\ell+1})
\Psi^{\mathfrak{gl}_{\ell+1}}_{\underline{\l}}
(\underline{p}_{\ell+1})=(\sum_{ I_r}\prod\limits_{i\in I_r}
q^{\l_i} )\,\, \Psi^{\mathfrak{gl}_{\ell+1}}_{\underline{\l}}
(\underline{p}_{\ell+1}),\ee where $\underline{\l}=(\lambda_1,
\ldots,\lambda_{\ell+1})$.

Let $\,{\cal P}^{(\ell+1)}\,$  be a set of collections of 
integers $\,p_{i,j}\in \ZZ\,,\,i=1,\ldots,\ell+1\,,\,j=1,\ldots,i\,$
satisfying the conditions $\,p_{i+1,j}\leq p_{i,j}\leq
p_{i+1,j+1}\,$ with fixed $\,p_{\ell+1,i}\,,\,i=1,\ldots,\ell+1\,$.
Thus $\,{\cal P}^{(\ell+1)}\,$  is a set of 
Gelfand-Zetlin patterns corresponding to an irreducible
finite-dimensional representation of $\,GL(\ell+1,\C)\,$ 
(see e.g. \cite{ZS}). 
We denote by $\,{\cal P}_{\ell+1,\ell}\subset {\cal P}^{(\ell+1)}\,
$ the subset $\,p_{\ell+1,i}\leq p_{\ell,i}\leq
p_{\ell+1,i+1}\,,\,i=1,\ldots,\ell$ of ${\cal P}^{(\ell+1)}\, $.

\begin{te}\label{mainTH}
The following  function is a solution  of the eigenfunction problem (\ref{eiglat})
 \be\label{main}
\Psi^{\mathfrak{gl}_{\ell+1}}_{\underline{\l}}
(\underline{p}_{\ell+1})\,=\, \sum_{p_{k,i}\in{\cal
P}^{(\ell+1)}}\,\,
\prod_{k=1}^{\ell+1}\,\, q^{\l_k(\sum_{i=1}^k p_{k,i}-\sum_{i=1}^{k-1} p_{k-1,i})}\\
\times\frac{\prod\limits_{k=2}^{\ell}\prod\limits_{i=1}^{k-1}
(p_{k,i+1}-p_{k,i})_q!}
{\prod\limits_{k=1}^{\ell}\prod\limits_{i=1}^k
(p_{k,i}-p_{k+1,i})_q!\,\, (p_{k+1,i+1}-p_{k,i})_q!}, \hspace{0.5cm}
p_{\ell+1,1}\leq \cdots \leq p_{\ell+1,\ell+1}, \\
\Psi^{\mathfrak{gl}_{\ell+1}}_{\underline{\l}}
(\underline{p}_{\ell+1})\,=\,0,\qquad \mbox{otherwise}. \ee Here we
set $(n)_q!=(1-q)...(1-q^n)$.

\end{te} The proof of the theorem  will be given in Section 4.

\begin{ex} Let $\mathfrak{g}=\mathfrak{gl}_{2}$ and
  $(p_{2,1},p_{2,2})\in\ZZ^{2}$.  The function
$$
\Psi_{\l_1,\l_2}^{{\mathfrak gl}_2}(p_{2,1},p_{2,2})
=\sum_{p_{2,1}\leq p_{1,1}\leq p_{2,2}}\frac{
q^{\l_1\,p_{1,1}} q^{\l_2(p_{2,1}+p_{2,2}-p_{1,1})}}
{(p_{1,1}-p_{2,1})_q!(p_{2,2}-p_{1,1})_q!},\qquad  p_{2,1}\leq p_{2,2}\,,
$$
$$
\Psi_{\l_1,\l_2}^{{\mathfrak gl}_2}(p_{2,1},p_{2,2})=0, \qquad p_{2,1}>p_{2,2}\,,
$$
is a  common eigenfunction of   commuting Hamiltonians
$$
{\cal H}_1^{{\mathfrak gl}_2}=
T_1+(1-q^{p_{2,2}-p_{2,1}+1})T_2,\qquad  {\cal H}_2^{{\mathfrak
    gl}_2}=T_1T_2.
$$
\end{ex}

Note that the  formula (\ref{main}) can be easily rewritten in the recursive
form.
\begin{cor} The following recursive relation holds
 \be\label{qtodarec}
\Psi^{\mathfrak{gl}_{\ell+1}}_{\underline{\l}}
(\underline{p}_{\ell+1})\,=\,\sum_{p_{\ell,i}\in{\cal
P}_{\ell+1,\ell}}\,
\Delta(\underline{p}_{\ell})\,q^{\l_{\ell+1}(\sum_{i=1}^{\ell+1}
p_{\ell+1,i}-\sum_{i=1}^\ell p_{\ell,i})}\,\,
\CQ_{\ell+1,\ell}(\underline{p}_{\ell+1},\underline{p}_{\ell}|q)
\Psi^{\mathfrak{gl}_{\ell}}_{\underline{\l}'}(\underline{p}_{\ell})\,,
\ee where
$$
\CQ_{\ell+1,\ell}(\underline{p}_{\ell+1},\underline{p}_{\ell}|q)\,=\,
\frac{1} {\prod\limits_{i=1}^{\ell} (p_{\ell,i}-p_{\ell+1,i})_q!\,\,
(p_{\ell+1,i+1}-p_{\ell,i})_q!}\,,
$$
and
$$
\Delta(\underline{p}_{\ell})=
\prod_{i=1}^{\ell-1}(p_{\ell,i+1}-p_{\ell,i})_q!\,,
$$
where  the notations $\underline{\l}=(\l_1,\ldots,\l_{\ell+1})$,
$\underline{\l}'=(\l_1,\ldots,\l_{\ell})$ are used.
\end{cor}

\begin{rem}
In the limit $q\to 1$ the $q$-deformed Toda chain eigenfunction equations
reduce to ordinary Toda chain eigenfunction equations. The solution
\eqref{main} in the limit $q\to 1$ reduces to an integral
representation of $\mathfrak{gl}_{\ell+1}$-Whittaker functions due
to Givental \cite{Gi} 
\bqa\label{giv11}
\Psi_{\underline{\lambda}}^{\mathfrak{gl}_{\ell+1}}
(x_1,\ldots,x_{\ell+1})=\,
 \int_{\RR^{\frac{\ell(\ell+1)}{2}}}
\prod_{k=1}^{\ell}\prod_{i=1}^kdx_{k,i}\,\,
e^{\mathcal{F}^{\mathfrak{gl}_{\ell+1}}(x) }, \eqa where
\bqa\label{intrep11} \mathcal{F}^{{\mathfrak{gl}}_{\ell+1}}(x)=
\imath\sum\limits_{k=1}^{\ell+1} \lambda_k\Big(\sum\limits_{i=1}^{k}
x_{k,i}-\sum\limits_{i=1}^{k-1}x_{k-1,i}\Big)-
\sum\limits_{k=1}^{\ell} \sum\limits_{i=1}^{k}
\Big(e^{x_{k,i}-x_{k+1,i}}+e^{x_{k+1,i+1}-x_{k,i}}\Big), \nonumber
\eqa  $\underline{\l}=(\l_1,\ldots,\l_{\ell+1})$ and
$x_i:=x_{\ell+1,i},\,\,\,i=1,\ldots,\ell+1$.
\end{rem}

\begin{prop}\label{qCS}
The common eigenfunction \eqref{main} of $q$-deformed Toda chain
allows the following representation for 
$p_{\ell+1,1}\leq p_{\ell+1,2}\leq\ldots p_{\ell+1,\ell+1}$
\be\label{intc}$$ \Psi^{\mathfrak{gl}_{\ell+1}}_{\underline{\l}}
(\underline{p}_{\ell+1})\,=\,\Tr_{V}\, q^{d}\prod_{i=1}^{\ell+1}
q^{\l_i\,E_{i,i}}, \ee where $V$ is a $\IC^*\times
GL(\ell+1,\IC)$-module, $E_{i,i}$, $i=1,\ldots \ell+1$ are Cartan
generators of $\mathfrak{gl}_{\ell+1}={\rm Lie}(GL(\ell+1,\IC))$ and
$d$ is a generator of ${\rm Lie}(\IC^*)$.
\end{prop}
{\it Proof}: It is useful to rewrite \eqref{main} in the following
form
$$
\Psi^{\mathfrak{gl}_{\ell+1}}_{\underline{\l}}
(\underline{p}_{\ell+1})\,=\,\Delta(\underline{p}_{\ell+1})^{-1}\,\,
\sum_{p_{k,i}\in\CP^{(\ell+1)}}\prod_{k=1}^{\ell}\,q^{\l_{k+1}(\sum_i
p_{k+1,i}-\sum_j p_{k,i})} \prod_{i=1}^k\,
{p_{k+1,\,i+1}-p_{k+1,\,i}\choose p_{k,\,i}-p_{k+1,\,i}}_{\!\!q}\,,
$$
$$
\Delta(\underline{p}_{\ell+1})\,=\,
\prod_{j=1}^{\ell}(p_{\ell+1,j+1}-p_{\ell+1,j})_q!\,,
$$
where 
$$
{n \choose k}_{\!\!q}\,=\,\frac{(n)_q!}{(n-k)_q!\,\,(k)_q!}.
$$
Now taking into account  the identities 
$$
{n \choose k}_{\!\!q}\,=\,
{n-1 \choose k-1}_{\!\!q}\,+\,q^k{n-1 \choose k}_{\!\!q},\qquad 
\frac{1}{(1-q^n)}=\sum_{k=0}^{\infty} q^{kn},
$$
one obtains a representation of $\Psi^{\mathfrak{gl}_{\ell+1}}_{\underline{\l}}
(\underline{p}_{\ell+1})$ as a sum of terms $q^{N}z_1^{m_1}\cdots
z_{\ell+1}^{m_{\ell+1}}$ with positive integer coefficients. Let
us note that $q$-deformed Toda chain eigenfunction equations 
depend on variables $z_i=q^{\l_i}$  only through characters of
fundamental representations of $\mathfrak{gl}_{\ell+1}$. Checking
that the initial conditions leading to the solution \eqref{main} can
be also expressed through characters and using the expansion with
positive integral coefficients discussed above, one obtains the
representation \eqref{intc} $\Box$

\begin{rem} The following representation holds 
for $p_{\ell+1,1}\leq p_{\ell+1,2}\leq\ldots p_{\ell+1,\ell+1}$
\be\label{inrep1}
\widetilde{\Psi}^{\mathfrak{gl}_{\ell+1}}_{\underline{\l}}
(\underline{p}_{\ell+1})\,=\,\Delta(\underline{p}_{\ell+1})\,\,
\Psi^{\mathfrak{gl}_{\ell+1}}_{\underline{\l}}
(\underline{p}_{\ell+1})=\Tr_{V_f}\,q^{d}\prod_{i=1}^{\ell+1}
q^{\l_i\,E_{i,i}}, \ee where $V_f$ is a {\it finite-dimensional}
$\IC^*\times GL(\ell+1,\IC)$-module (see Proposition 3.4 for 
explicite description of $V_f$). The module $V$ entering \eqref{intc} and the module
$V_f$ entering \eqref{inrep1} have a more refined structure under
the action of (quantum)  affine Lie algebras and will be discussed
in an other part of this paper \cite{GLO2}.
\end{rem}

\section{Non-standard  limits}

Besides the limit $q\to 1$ recovering classical
$\mathfrak{gl}_{\ell+1}$-Whittaker function as a solution of
$\mathfrak{gl}_{\ell+1}$-Toda chain there are other interesting
limits elucidating the meaning of $q$-deformed Toda chain equations.
In the limit $q\to 0$  the $q$-deformed
$\mathfrak{gl}_{\ell+1}$-Whittaker functions are given by characters
of irreducible representations of $\mathfrak{gl}_{\ell+1}$. This
allows to identify the Whittaker functions with $p$-adic Whittaker
functions according to Shintani-Casselman-Shalika formula \cite{Sh},
\cite{CS}. There is also a non-standard $q\to 1$ limit  which
clarifies the recursive structure of $q$-deformed
$\mathfrak{gl}_{\ell+1}$-Whittaker functions.

\subsection{The limit $q\to 0$}

In this subsection we discuss a limit $q\to 0$ of the constructed
$q$-deformed Whittaker functions (we restrict Whittaker function to
the domain $\{p_{\ell+1,1}\leq\ldots\leq p_{\ell+1,\ell+1}\}$
where it is non-trivial). We will show that in the domain 
$\{p_{\ell+1,1}\leq\ldots\leq p_{\ell+1,\ell+1}\}$ the system of
equations for a common eigenfunctions of  $q$-deformed Toda chain
Hamiltonians reduces to Pieri formulas (particular case of
Littlewood-Richardson rules)   for the decomposition of the tensor
product of an arbitrary finite-dimensional representation and a
fundamental representation of $\mathfrak{gl}_{\ell+1}$.

\begin{prop}\label{qzeropr}

1. In the limit $q\to 0$, the solution \eqref{main} is given in the
domain $p_{\ell+1,1}\leq\ldots\leq p_{\ell+1,\ell+1}$  by
\be\label{qzero} \Psi^{\mathfrak{gl}_{\ell+1}}_{\underline{\l}}
(\underline{p}_{\ell+1})|_{q\to 0}:=
\chi^{\mathfrak{gl}_{\ell+1}}_{\underline{p}_{\ell+1}}
(\underline{z})\,=\,\sum_{p_{k,i}\in{\cal P}^{\ell+1}}\,\,
\prod_{k=1}^{\ell+1} z_k^{(\sum_{i=1}^k p_{k,i}-\sum_{i=1}^{k-1}
p_{k-1,i})}\,,\,\ee where we set
$z_i=q^{\l_i}\,,\,i=1,\ldots,\ell+1\,.$ 

2. Functions $\chi_{\underline{p}_{\ell+1}}(z)$ satisfy the
following set of difference equations \be\label{Htozero}
\chi^{\mathfrak{gl}_{\ell+1}}_{r}(\underline{z})\,
\chi^{\mathfrak{gl}_{\ell+1}}_{\underline{p}_{\ell+1}}
(\underline{z})\,=\,\sum_{I_r}\,\,
\chi^{\mathfrak{gl}_{\ell+1}}_{\underline{p}_{\ell+1}+I_r}
(\underline{z})\,, \ee where 
$\underline{z}=(z_1,z_2,\ldots,z_{\ell+1})$.

3. The functions
$\chi^{\mathfrak{gl}_{\ell+1}}_{\underline{p}_{\ell+1}}(\underline{z})$
can be identified with characters of irreducible
finite-dimensional representations of $GL_{\ell+1}$  corresponding
to partitions $p_{\ell+1,1}\leq\ldots\leq p_{\ell+1,\ell+1}$.
\end{prop}

{\it Proof}: The relations \eqref{qzero} and \eqref{Htozero} follows
directly from the similar relations for generic $q$. To prove the
last statement note that \eqref{qzero} can be identified  with
expression for characters of irreducible finite-dimensional
representations of $GL_{\ell+1}$ obtained using the Gelfand-Zetlin
bases (see e.g. \cite{ZS}).
 Let $\{p_{ij}\}$, $i=1,\ldots,\ell+1$, $j=1,\ldots, i$   be a  Gelfand-Zetlin (GZ)
pattern ${\cal P}^{(\ell+1)}$ i.e.  satisfy the conditions
$p_{i+1,j}\leq p_{i,j}\leq p_{i+1,j+1}$. Irreducible
finite-dimensional representation can be realized in a vector space
with the bases $v_{\underline{p}}$ parametrized by GZ patterns
$\{p_{ij}\}$ with fixed $p_{\ell+1,i}$. Action of Cartan generators
on $v_{\underline{p}}$ is given by \be z_1^{E_{11}}\,z_2^{E_{22}}
\cdots z_{\ell+1}^{E_{\ell+1,\ell+1}}\,v_{\underline{p}}=
z_1^{s_1}\,z_2^{s_2-s_1}\,\cdots
z_{\ell+1}^{s_{\ell+1}-s_{\ell}}\,v_{\underline{p}},\qquad
s_k=\sum_{i=1}^k\,p_{ki}. \ee Thus we have for the character
\be\label{character} \chi^{\mathfrak{gl}_{\ell+1}}
_{p_{\ell+1,1},\ldots ,p_{\ell+1,\ell+1}}(z_1,\ldots ,z_{\ell+1})=
\sum_{ p_{k,i}\in{\cal P}^{(\ell+1)}}\,\prod_{k=1}^{\ell+1}
z_k^{(\sum_{i=1}^k p_{k,i}-\sum_{i=1}^{k-1} p_{k-1,i})}\, \ee
$\Box$

\begin{rem}
The second identity in Proposition \ref{qzeropr}  is known as the  Pieri
formula  (see e.g. \cite{FH}, Appendix A). Thus the $q$-deformed
Toda chain equations can be considered as $q$-deformations of Pieri's formula.
There is a generalization of $q$-Toda chain relations providing
a  $q$-version of a general Littlewood-Richardson rule.
\end{rem}

The expressions for the characters in GZ representation have  an
obvious recursive structure which is a $q\to 0$ limit of
(\ref{qtodarec}).

\begin{cor}
Characters satisfy the following recursive relation \be
\chi^{\mathfrak{gl}_{\ell+1}}_{p_{\ell+1,1},\ldots
,p_{\ell+1,\ell+1}} (z_1,\ldots ,z_{\ell+1})=
\sum_{p_{\ell,i}\in{\cal P}_{\ell+1,\ell}}
   \,\,\, z_{\ell+1}^{ \sum_{i=1}^{\ell+1}
  p_{\ell+1,i}-\sum_{i=1}^{\ell}  p_{\ell,i}}
 \chi^{\mathfrak{gl}_{\ell}}_{p_{\ell,1},\ldots ,p_{\ell,\ell}}
 (z_1,\ldots ,z_{\ell}),
\ee where sum goes over $\underline{p}_{\ell}=(p_{\ell,1},\ldots
,p_{\ell,\ell})$ satisfying the GZ conditions:
$$
p_{\ell+1,i}\leq p_{\ell,i}\leq p_{\ell+1,i+1}.
$$
\end{cor}
Note that these recursive relations can be derived using the
classical Cauchy-Littlewood formula
\bqa\label{CauchyA} C_{\ell+1,m+1}(x,y)=
\prod_{i=1}^{\ell+1}\prod_{j=1}^{m+1}\,\frac{1}{1-x_i y_j}=
\sum_{\Lambda}\,\chi^{\mathfrak{gl}_{\ell+1}}_{\Lambda}(x)\,\,
\chi^{\mathfrak{gl}_{m+1}}_{\Lambda}(y) ,\eqa where the sum goes
over Young diagrams $\Lambda$ of $\mathfrak{gl}_{m+1}$ and
$\chi^{\mathfrak{gl}_{\ell+1}}_{\Lambda}(x)=
\chi^{\mathfrak{gl}_{\ell+1}}_{\Lambda} (x_1,\ldots , x_{\ell+1})$
are characters of irreducible finite-dimensional representation of
$GL(\ell+1,\IC)$ corresponding to  Young diagrams $\Lambda$.

\begin{prop} The following integral relations for the characters
$\chi^{\mathfrak{gl}_{\ell+1}}_{\Lambda}(x)$ hold
\be\label{reccharGZ} \chi^{\mathfrak{gl}_{\ell+1}}_{\Lambda}(x)\,=\,
\oint_{y_1=0}\cdots \oint_{y_{\ell}=0}\, \prod_{i=1}^{\ell}\,\frac{dy_i}{2\pi\i y_i} \,
C_{\ell+1,\ell}(x,y^{-1})
\chi^{\mathfrak{gl}_{\ell}}_{\Lambda}(y)\,\Delta(y), \ee
$$
\chi^{\mathfrak{gl}_{\ell+1}}_{\Lambda}(x)\,=\,
\oint_{y_1=0}\cdots \oint_{y_{\ell+1}=0}\,
\prod_{i=1}^{\ell+1}\,\frac{dy_i}{2\pi\i y_i} \,
C_{\ell+1,\ell+1}(x,y^{-1})
\chi^{\mathfrak{gl}_{\ell+1}}_{\Lambda}(y)\,\Delta(y),
$$
$$\chi^{\mathfrak{gl}_{\ell+1}}_{\Lambda+(\ell+1)^k}(x)=\big(\prod_{j=1}^{\ell+1}x_j^k\big)
\chi^{\mathfrak{gl}_{\ell+1}}_{\Lambda}(x).$$
\end{prop}
The relations above  allow to obtain
 character of irreducible finite-dimensional representation of
$GL(\ell+1,\IC)$ corresponding to any Young diagram $\Lambda$.

\begin{rem} The relations above can be obtained from similar
  relations for Macdonald polynomials in the  limit $t\to 0$, $q\to
  0$. These recursion relations are analogs of Mellin-Barnes recursion
relations for classical Whittaker functions (see
\cite{KL1}, \cite{GKL}, \cite{GLO} for details).
\end{rem}

According to  Shintani-Casselman-Shalika formula, the $p$-adic
Whittaker function corresponding to a Lie group $G$ is equal to the
character of the  Langlands dual Lie group  ${}^LG$ acting in
an irreducible finite-dimensional representation \cite{Sh},
\cite{CS}. Thus according to Proposition \ref{qzero} we  can
consider $\mathfrak{gl}_{\ell+1}$-Whittaker functions at $q\to 0$ as
an incarnation of $p$-adic Whittaker functions (this is in complete
agreement with the results of \cite{GLO3}.
Moreover, taking into  account Proposition \ref{qCS} one can
consider the main result of this paper as a generalization of
Shintani-Casselman-Shalika formula to a $q$-deformed case including
a limiting case of classical $\mathfrak{gl}_{\ell+1}$-Whittaker
functions. This interpretation of classical Whittaker functions
evidently deserves further attention.

\subsection{Modified  limit $q\to 1$}

In this subsection we consider a modified limit $q\to 1$ leading to
a very simple degeneration of $q$-deformed Toda chain. In this limit
$q$-deformed Toda chain can be easily solved. Moreover  the form of
the solution makes the recursive expressions \eqref{qtodarec} for
$q$-deformed Toda chain solution very natural.

Let us redefine the  $q$-deformed Toda chain Hamiltonians and
their common eigenfunctions as follows (we assume 
$p_{\ell+1,1}\leq\ldots \leq p_{\ell+1,\ell+1}$ below )
$$
\CJ^{\mathfrak{gl}_{\ell+1}}_r=
\Delta(\underline{p}_{\ell+1})\,\,\,{\ch}_r^{\mathfrak{gl}_{\ell+1}}\,\,\,
\Delta(\underline{p}_{\ell+1})^{-1}\,,
$$
$$
\widetilde{\Psi}_{\underline{\lambda}}^{\mathfrak{gl}_{\ell+1}}
(\underline{p}_{\ell+1})\,=\, \Delta(\underline{p}_{\ell+1})\,\cdot
\Psi^{\mathfrak{gl}_{\ell+1}}_{\underline{\lambda}}(\underline{p}_{\ell+1}),\qquad
\Delta(\underline{p}_{\ell+1})\,=\,
\prod_{j=1}^{\ell}(p_{\ell+1,j+1}-p_{\ell+1,j})_q!\,.
$$
Explicitly we have \be\label{Modifqto1Hamiltonians}
\CJ_r^{\mathfrak{gl}_{\ell+1}}\,=\,\sum_{I_r}\,\bigl(
\widetilde{X}_{i_1}^{1-\delta_{i_2-i_1,\,1}}\cdot\ldots\cdot
\widetilde{X}_{i_{r-1}}^{1-\delta_{i_r-i_{r-1},\,1}}\cdot
\widetilde{X}_{i_r}^{1-\delta_{i_{r+1}-i_r,\,1}}\bigr)
T_{i_1}\cdot\ldots\cdot T_{i_r},
\ee
 and we assume $i_{r+1}=\ell+2$
and $ \,\widetilde{X}_i\,=\,1-q^{p_{i+1}-p_i}\, $.
Let us now take the limit $q\to 1$
$$
\widetilde{\psi}^{(\ell+1)}_{\underline{\lambda}}(\underline{p}_{\ell+1})\,=\,
\lim_{q\to 1}\,
\widetilde{\Psi}^{\mathfrak{gl}_{\ell+1}}_{\underline{\lambda}}
(\underline{p_{\ell+1}}),\hspace{1.5cm}
h_r^{(\ell+1)}=\lim_{q\to 1}
\,\,\CJ_r^{{\mathfrak gl}_{\ell+1}}\, . $$
 We have  $\lim_{q\to1}(1-q^n)=0$, and therefore we  obtain
from (\ref{comm})
$$
h_r^{(\ell+1)}\,=\,T_{\ell+2-r}\cdot\ldots\cdot T_{\ell+1}.
$$
Now the eigenfunction problem is easily solved.

\begin{prop}\label{qzerochar} 1. The function
\be\label{eigchar} \widetilde{\psi}^{(\ell+1)}_{\underline{\lambda}}
(\underline{p}_{\ell+1})\,=\,
\Big(\chi_{\ell+1}^{\mathfrak{gl}_{\ell+1}}
(\underline{z})\Big)^{p_{\ell+1,\,1}} \prod_{i=1}^\ell\,
\Big(\chi_{\ell+1-i}^{\mathfrak{gl}_{\ell+1}}
(\underline{z})\Big)^{p_{\ell+1,\,i+1}-p_{\ell+1,\,i}} \ee
 is an eigenfunction of the family of mutually commuting difference
operators
 \be\label{eigq1}
h_r^{(\ell+1)}
\widetilde{\psi}^{(\ell+1)}_{\underline{\lambda}}(\underline{p}_{\ell+1})\,=\,
\chi^{\mathfrak{gl}_{\ell+1}}_{r}(\underline{z})\,
\widetilde{\psi}^{(\ell+1)}_{\underline{\lambda}}(\underline{p}_{\ell+1}),
\ee where $\chi_{r}^{\mathfrak{gl}_{\ell+1}}(\underline{z})$ is the
character of fundamental representation
$V_{\omega_r}=\bigwedge^r\mathbb{C}^{\ell+1}$:
$$\chi^{\mathfrak{gl}_{\ell+1}}_{r}(\underline{z})\,=\,
\sum_{I_r}\,z_{i_1}\cdots
z_{i_r}\,,\,\,\,z_i=q^{\lambda_i}\,,\,\,\,\,\,i=1,\ldots,\ell+1,
$$
and
$$ h_r^{(\ell+1)}\,=\,T_{\ell+2-r}\cdot\ldots\cdot
T_{\ell+1}\,,\,\,\,r=1,\ldots,\ell+1.
$$

2. In the domain $ p_{\ell+1,1}\leq\ldots\leq p_{\ell+1,\ell+1}$ the
following recursive relation holds  \be\label{eigq11}
\widetilde{\psi}^{(\ell+1)}_{\underline{\lambda}}(\underline{p}_{\ell+1})\,=\,
\sum_{ p_{\ell,i}\in {\cal P}_{\ell+1,\ell}}\,
z_{\ell+1}^{\sum_{i=1}^{\ell+1} p_{\ell+1,i}-\sum_{i=1}^{\ell}
p_{\ell,i}}
\prod_{i=1}^\ell\,{p_{\ell+1,\,i+1}-p_{\ell+1,\,i}\choose
p_{\ell,\,i}-p_{\ell+1,\,i}}_{\!}\,\cdot
\widetilde{\psi}^{(\ell)}_{\underline{\lambda}'}(\underline{p}_\ell),
\ee where $\underline{\lambda}'=(\lambda_1,\ldots,\lambda_{\ell})$
and $z_{\ell+1}=q^{\lambda_{\ell+1}}.$
\end{prop}

\noindent{\it Proof}: The identity \eqref{eigq1} follows from the construction. 
Let us prove that \eqref{eigq11}  follows from  \eqref{eigq1}. Denote
$(z_1,\ldots,z_{\ell+1})=(q^{\l_1},\ldots,q^{\l_{\ell+1}})$. 
Using the relation 
$$
\chi^{\mathfrak{gl}_{\ell+1}}_r(\underline{z})\,=\,
\chi^{\mathfrak{gl}_\ell}_r(\underline{z}')+
z_{\ell+1}\cdot\chi^{\mathfrak{gl}_\ell}_{r-1}(\underline{z}'),\qquad 
r\,=\, 1,\ldots,\ell+1,
$$
we have
\be\Big(\chi_{\ell+1}^{\mathfrak{gl}_{\ell+1}}
(\underline{z})\Big)^{p_{\ell+1,\,1}}
\prod_{i=1}^{\ell}\,\Big(\chi_{\ell+1-i}^{\mathfrak{gl}_{\ell+1}}
(\underline{z})\Big)^{p_{\ell+1,\,i+1}-p_{\ell+1,\,i}}\\=
\Big(z_{\ell+1}\chi_{\ell}^{\mathfrak{gl}_\ell}
(\underline{z}')\Big)^{p_{\ell+1,1}}
\prod_{i=1}^{\ell}\,\,\,\sum_{p_{\ell,i}=p_{\ell+1,i}}^{p_{\ell+1,i+1}}
\Big(z_{\ell+1}\chi_{\ell-i}^{\mathfrak{gl}_\ell}
(\underline{z}')\Big)^{p_{\ell+1,i+1}-p_{\ell,i}}\\ \cdot
\Big(\chi_{\ell+1-i}^{\mathfrak{gl}_\ell}
(\underline{z}')\Big)^{p_{\ell,i}-p_{\ell+1,i}}
\cdot{p_{\ell+1,i+1}-p_{\ell+1,\,i}\choose
p_{\ell,\,i}-p_{\ell+1,\,i}}\\= \sum_{p_{\ell,i}\in {\cal
P}_{\ell+1,\ell}}z_{\ell+1}^{\sum_i p_{\ell+1,i+1}-\sum_i
p_{\ell,i}}\prod_{i=1}^{\ell}{p_{\ell+1,i+1}-p_{\ell+1,\,i}\choose
p_{\ell,\,i}-p_{\ell+1,\,i}}\\ \cdot
\Big(\chi^{\mathfrak{gl}_\ell}_{\ell}
(\underline{z}')\Big)^{p_{\ell,\,1}}
\prod_{i=1}^{\ell-1}\,\Big(\chi_{\ell-i}^{\mathfrak{gl}_\ell}
(\underline{z}')\Big)^{p_{\ell,\,i+1}-p_{\ell,\,i}}\\=
\sum_{p_{\ell,i}\in {\cal P}_{\ell+1,\ell}} z_{\ell+1}^{\sum_i
p_{\ell+1,i+1}-\sum_i
p_{\ell,i}}\prod_{i=1}^{\ell}{p_{\ell+1,i+1}-p_{\ell+1,\,i}\choose
p_{\ell,\,i}-p_{\ell+1,\,i}}\,\,
\widetilde{\psi}_{\underline{\lambda}'}^{(\ell)}(\underline{p}_{\ell}).
\ee Thus we obtain the recursion formula described in the
proposition $\Box$

\begin{rem} The functions
$$
\psi^{(\ell+1)}_{\underline{\lambda}}(\underline{p}_{\ell+1})\,=\,
\Delta^{-1}(\underline{p}_{\ell+1})
\widetilde{\psi}^{(\ell+1)}_{\underline{\lambda}}(\underline{p}_{\ell+1}),
$$
satisfy following recursive relations
 \be\label{simplerec}
{\psi}^{(\ell+1)}_{\underline{\l}}(\underline{p}_{\ell+1})\,=\,
\sum_{p_{\ell,i}\in{\cal P}_{\ell+1,\ell}}\, z_{\ell+1}^{\sum_i
p_{\ell+1,i}-\sum_i p_{\ell,i}}
\prod_{i=1}^{\ell-1}\frac{(p_{\ell,i+1}-p_{\ell,i})!}
{(p_{\ell,i}-p_{\ell+1,i})!\,\, (p_{\ell+1,i+1}-p_{\ell,i})!}\,
{\psi}^{(\ell)}_{\underline{\lambda}'}(\underline{p}_{\ell}). \ee
This makes the formula \eqref{qtodarec} for the solution of $q$-deformed
Toda chain slightly less mysterious.
\end{rem}

\begin{prop} The following representation holds
\be\label{inrep2}
\widetilde{\psi}^{(\ell+1)}_{\underline{\l}}(\underline{p}_{\ell+1})
\,=\,\Tr_{V_f}\,\prod_{i=1}^{\ell+1} q^{\l_i\,E_{i,i}}, \ee where
\be\label{decomp} V_f\,=\,V_{\omega_1}^{\otimes
(p_{\ell+1,\ell+1}-p_{\ell+1,\ell})}\otimes \cdots \otimes
V_{\omega_{\ell}}^{\otimes (p_{\ell+1,2}-p_{\ell+1,1})} \otimes
V_{\omega_{\ell+1}}^{\otimes p_{\ell+1,1}}, \ee where
$V_{\omega_n}=\wedge^{n}\IC$ are the fundamental representations of
$GL(\ell+1,\IC)$.
\end{prop}
\noindent {\it Proof}: Obvious consequence of the Proposition
\ref{qzerochar} $\Box$

The module $V_f$ entering \eqref{inrep1} is isomorphic to
\eqref{decomp} as $GL(\ell+1,\IC)$-module but has a more refined
structure under the action of quantum affine Lie algebras 
and will be discussed in  \cite{GLO2}.

\section{ Proof of Theorem \ref{mainTH}}

In this section we provide a proof of Theorem \ref{mainTH}.
To derive explicit expression \eqref{main} for $q$-deformed
 $\mathfrak{gl}_{\ell+1}$-Whittaker function we take a  limit $t\to
 \infty$  of recursive relations (\ref{macrec})
 for Macdonald polynomials.
We start with some useful relations
that will be used in the proof of Theorem \ref{mainTH}.

\subsection{$q$-deformed Toda chain from Macdonald-Ruijsenaars system}

In this subsection we demonstrate that
quantum Hamiltonians of $q$-deformed Toda chain arise as a limit of
Macdonald-Ruijsenaars operators when  $t\to \infty$. Let us take
$t=q^{-k}$.

\begin{prop}\label{escalc}
 The following relations hold
 $$
\CH_r^{\mathfrak{gl}_{\ell+1}}\,=\,\lim_{k\to \infty}
H_{r,k}^{\mathfrak{gl}_{\ell+1}}\,=\\
\sum_{I_r}\,\bigl(X_{i_1}^{1-\delta_{i_1,\,1}}\cdot
X_{i_2}^{1-\delta_{i_2-i_1,\,1}}\cdot\ldots\cdot
X_{i_r}^{1-\delta_{i_r-i_{r-1},\,1}}\bigr) T_{i_1}\cdot\ldots\cdot
T_{i_r},
$$
where
$$
H_{r,k}^{\mathfrak{gl}_{\ell+1}}\,=\,
D(x)^{-1}\,H_r^{\mathfrak{gl}_{\ell+1}}(x_iq^{k\,i})\,D(x),\qquad
D(x) =\prod_{i=1}^{\ell+1}x_i^{-k(\ell+1-i)},
$$
and the sum is over subsets
$I_r=\{i_1<i_2<\ldots<i_r\}\subset\{1,2,\ldots,\,\ell+1\}$. We take
$\,X_i=1-x_ix_{i-1}^{-1}\,$,$\,i=2,\ldots,\ell+1\,$, with
$\,X_1=1\,$ and $\,T_i\,x_j=q^{\delta_{i,j}}\,x_j\,T_i$.
\end{prop}

\noindent {\it Proof}: Make a change of variables $x_i$: $x_i\longmapsto
x_it^{-i}$, $i=1,\ldots,\ell+1$. Then for any $i$ and any $I_r$,
containing $i$ we have:
\be\label{eterm}
\Big(\prod_{j\notin I_r}
\frac{tx_i-x_j}{x_i-x_j}\Big)\longmapsto
\Big(t^{b_{r,i}}\prod_{j>i}\frac{x_i-x_jt^{i-1-j}}{x_i-x_jt^{i-j}}\times
\frac{x_i-x_{i-1}}{x_i t^{-1}-x_{i-1}}\times
\prod_{j<i-1}\frac{x_it^{j+1-i}-x_j}{x_it^{j-i}-x_j}\Big),
\ee where
$b_{r,i}=|\{j\notin I_r|j>i\}|$.
Making  a substitution $t=q^{-k}$ and conjugating the Hamiltonians
$H_r^{\mathfrak{gl}_{\ell+1}}$ by
$$
D(x)=\prod_{i=1}^{\ell+1}x_i^{-k\varrho_i},
$$
leads to a  multiplication
of each term \eqref{eterm} in the sum
(\ref{macdif}) by $\prod_{i\in I_r} q^{-k\varrho_i}$,
$\varrho_i:=\ell+1-i$.  Taking into account that for any $i$ and for any
subset $I_r$ containing $i$ one has $\sum_{i\in
I_r}\bigl(\varrho_i-b_{r,i}\bigr)\,=\,r(r-1)/2$
 we obtain in the limit
$k\to\infty$
$$
\CH_r^{\mathfrak{gl}_{\ell+1}}\,=\,\sum_{I_r}\,
\bigl(X_{i_1}^{1-\delta_{i_1,\,1}}\cdot
X_{i_2}^{1-\delta_{i_2-i_1,\,1}}\cdot\ldots\cdot
X_{i_r}^{1-\delta_{i_r-i_{r-1},\,1}}\bigr) T_{i_1}\cdot\ldots\cdot
T_{i_r},
$$
 where $ X_i=1-x_ix_{i-1}^{-1},\,\,\,\,i=2,\ldots,\ell+1$
and $X_1=1$  $\Box$

\subsection{ The  Cauchy-Littlewood  kernel $C_{\ell+1,\ell}(x,y|q,t)$
for $q$-deformed Whittaker functions}

In this subsection by taking an appropriate limit
of the Cauchy-Littlewood  kernel for Macdonald polynomials
we derive its analog for $q$-deformed Whittaker functions
and verify the intertwining relations with $q$-deformed
Toda chain Hamiltonians.

Let $t=q^{-k}$, $\varrho_i=\ell+1-i$. Given the  Cauchy-Littlewood
kernel $C_{\ell+1,\ell}(x,y|q,t)$ (\ref{mackern}), define
new kernel as
 \be\label{baxkern} Q_{\ell+1,\ell}(x,y|q)=\lim_{k\to\infty}\,\Big\{
\prod_{i=1}^\ell(x_iy_{\ell+1-i})^{-k\varrho_i}\cdot R_k(q)\cdot
C_{\ell+1,\,\ell}(x,y|q,q^{-k})\Big\}, \ee where
$$
R_k(q)=\prod_{j=1}^k\frac{-q^{a_j}}{(1-q^j)^{2\ell}}\,\,,\hspace{1.5cm}
a_j=\frac{\ell(\ell+1)}{2}\Big(j+\frac{\ell-1}{3}k\,\Big).
$$

\begin{prop} Then the following explicit expression for
$Q_{\ell+1,\ell}(x,y|q)$ defined by \eqref{baxkern} holds:
 \be\label{kerinf}
Q_{\ell+1,\ell}(x,y|q)\,=\,\prod_{i=1}^\ell\prod_{n=1}^\infty
\frac{1-(x_iy_i)^{-1}q^n}{1-q^n}\cdot\prod_{i=1}^\ell\prod_{n=1}^\infty
\frac{1-x_{i+1}y_iq^{-1}q^n}{1-q^n}\,. \ee
\end{prop}

\noindent {\it Proof:}  Making substitution $x_i\to x_it^{-i}$,
$y_i\to y_it^i$ in $C_{\ell+1,\,\ell}(x,y|q,t)$  and taking $t=q^{-k}$
we have
$$
C_{\ell+1,\,\ell}(x,y|q,t)\,=\,\prod_{n=0}^\infty\prod_{i=1}^\ell
\frac{1-x_iy_iq^{n-k}}{1-x_iy_iq^n}
\frac{1-x_{i+1}y_iq^n}{1-x_{i+1}y_iq^{n+k}}
$$
$$
\times\prod_{i=3}^{\ell+1}\prod_{j=1}^{i-2}
\frac{1-x_iy_jq^{n+(i-j-1)k}}{1-x_iy_jq^{n+(i-j)k}}
\frac{1-x_{\ell+2-i}y_{\ell+1-j}q^{n+(j-i)k}}
{1-x_{\ell+2-i}y_{\ell+1-j}q^{n+(j+1-i)k}}\,.
$$
One encounters four types of factors which can be rewritten as
$$
\prod_{n=0}\frac{1-xyq^{n-k}}{1-xyq^n}=\prod_{j=1}^k(1-xyq^{-j})\,=\,
(xy)^k\prod_{j=1}^k(-q^{-j})\bigl(1-(xy)^{-1}q^j\bigr),
$$
$$
\prod_{n=0}\frac{1-xyq^n}{1-xyq^{n+k}}\,=\,
\prod_{n=0}^{k-1}(1-xyq^n)\,=\,
\prod_{j=1}^k\bigl(1-xyq^{-1}q^j\bigr),
$$
$$
\prod_{n=0}\frac{1-xyq^{n-(m+1)k}}{1-xyq^{n-mk}}\,=\,
\prod_{j=k+1}^{2k}(1-xyq^{-j})\,=\,
(xy)^k\prod_{j=k+1}^{2k}(-q^{-j})\bigl(1-(xy)^{-1}q^j\bigr),
$$
$$
\prod_{n=0}\frac{1-xyq^{n+mk}}{1-xyq^{n+(m+1)k}}\,=\,
\prod_{j=k+1}^{2k}\bigl(1-xyq^{-1}q^j\bigr).
$$
Now it is easy to take the limit $k\to \infty$ and obtain
\eqref{kerinf} $\Box$

Let us introduce a set of slightly modified mutually commuting
Hamiltonians \be\label{hadjoint}
\widetilde{\CH}_r^{\mathfrak{gl}_{\ell}}(y)\,=\,
\sum_{I_r}\,\bigl(Y_{i_1}^{1-\delta_{i_2-i_1,\,1}}\cdot\ldots\cdot
Y_{i_{r-1}}^{1-\delta_{i_r-i_{r-1},\,1}}\cdot
Y_{i_r}^{1-\delta_{i_{r+1}-i_r,\,1}}\bigr) T_{i_1}\cdot\ldots\cdot
T_{i_r},\ee where
$$
Y_i(y)\,=\,1-y_i y_{i+1}^{-1},\qquad 1\leq i<\ell,
\hspace{1.5cm}Y_{\ell}=1.
$$
We assume here $I_r=(i_1<i_2<\ldots<i_r)\subset\{1,2,\ldots,\ell\}$
and we set $i_{r+1}=\ell+1$.

\begin{prop}\label{intqxy}
The following  intertwining relations hold \be\label{limbraid}
{{\CH}}_k^{\frak{gl}_{\ell+1}}(x)\,Q_{\ell+1,\ell}(x,y|q)\,=\,
\Big\{\widetilde{\CH}^{\frak{gl}_{\ell}}_{k-1}(y)+
\widetilde{\CH}^{\frak{gl}_\ell}_k(y)\Big\}\,Q_{\ell+1,\ell}(x,y|q),\ee where
$k=1,\ldots,\ell+1$.  Here $\,Q_{\ell+1,\ell}(x,y|q)\,$,
$\,\CH_k^{\mathfrak{gl}_{\ell+1}}(x)\,$ and
$\,\,\widetilde{\CH}_k^{\mathfrak{gl}_{\ell}}(y)$ are defined by (\ref{kerinf}),
(\ref{comm}) and (\ref{hadjoint}) respectively.
\end{prop}
{\it Proof}: Direct calculation similar to the one used in the proof
of Proposition \ref{escalc} $\Box$

Let us introduce a function
$\CQ_{\ell+1,\ell}(\underline{p}_{\ell+1},\underline{p}_{\ell}|q)$
on the lattice $\IZ^{\ell+1}\times \IZ^{\ell}$ as follows
$$
\CQ_{\ell+1,\ell}(\underline{p}_{\ell+1},\underline{p}_{\ell}|q)=
Q_{\ell+1,\ell}(q^{p_{\ell+1,i}+i-1},q^{-p_{\ell,i}-i+1}|q).
$$

\begin{cor}
The following explicit expression for 
$\CQ_{\ell+1,\ell}(\underline{p}_{\ell+1},\underline{p}_{\ell}|q)$
holds
$$
\CQ_{\ell+1,\ell}(\underline{p}_{\ell+1},\underline{p}_{\ell}|q)=\frac{
\prod_{i=1}^{\ell}\Theta(p_{\ell,i}-p_{\ell+1,i})
\Theta(p_{\ell+1,i+1}-p_{\ell,i})} {\prod_{i=1}^{\ell}
(p_{\ell,i}-p_{\ell+1,i})_q!\,\, (p_{\ell+1,i+1}-p_{\ell,i})_q!},
$$
where $\Theta(n)=1$ when $n\geq 0$ and $\Theta(n)=0$ otherwise.
\end{cor}

\begin{prop}\label{interlat}
For any $k=1,\ldots,\ell+1$ the following  intertwining relations
hold \be\label{braiding}
{{\CH}}_k^{\frak{gl}_{\ell+1}}(\underline{p}_{\ell+1})\CQ_{\ell+1,\ell}
(\underline{p}_{\ell+1},\underline{p}_{\ell}|q)\,=\,
\Big\{\widetilde{\CH}^{\frak{gl}_{\ell}}_{k-1}(-\underline{p}_{\ell})+
\widetilde{\CH}^{\frak{gl}_{\ell}}_k (-\underline{p}_{\ell})\Big\}
\CQ_{\ell+1,\ell}(\underline{p}_{\ell+1},\underline{p}_{\ell}|q).\ee
\end{prop}
{\it Proof}: Follows from Proposition \ref{intqxy} $\,\Box$

\subsection{Pairing}

Define a pairing: \be\label{todpair} \<f,g\>_q=\oint_{y_1=0}\cdots
\oint_{y_{\ell}=0} \, \prod\limits_{i=1}^{\ell}\frac{dy_i}{2\pi\i
y_i} \Delta (y)f(y^{-1})g(y), \ee where
\be\label{todmeas}\Delta(y)\,=\,
\prod_{n=1}^\infty\prod_{i=1}^{\ell-1}
\frac{1-q^n}{1-y_{i+1}y_{i}^{-1}q^{n-1}}\,,\hspace{1cm}
f(y^{-1})\,:=\, f(y_1^{-1},\ldots,y_{\ell}^{-1}).\ee
\begin{prop}\label{selfadj}  Hamiltonians
$\CH_r^{\mathfrak{gl}_{\ell}}(y)$ and
$\widetilde{\CH}_r^{\mathfrak{gl}_\ell}(y)$ are adjoint with respect to the
pairing  (\ref{todpair})
$$
\<f,\ch_k^{\mathfrak{gl}_\ell} g\>_q\,=\,
\<\widetilde{\CH}_k^{\mathfrak{gl}_\ell} f,g\>_q,\qquad k=1,\ldots,\ell.
$$
\end{prop}

\noindent {\it Proof}: Let us adopt the following notations
$$
\ch_r^{\mathfrak{gl}_\ell}(y)\,=\,\sum_{I_r}
A_{I_r}(y)\,T_{I_r},\hspace{1.5cm}
\widetilde{\CH}_r^{\mathfrak{gl}_\ell}(y)\,=\,\sum_{I_r} B_{I_r}(y)\,T_{I_r},
$$ where
$T_{I_r}:=T_{i_1}T_{i_2}\cdot\ldots\cdot T_{i_r}$. One should prove
\be \oint_{y_1=0}\cdots \oint_{y_{\ell}=0}
\prod_{i=1}^{\ell}\frac{dy_i}{2\pi\i y_i}\Delta(y)
f(y^{-1})\sum_{I_r}T_{I_r} \cdot T_{I_r}^{-1} A_{I_r}(y)
T_{I_r}g(y)\\= \nonumber\oint_{y_1=0}\cdots \oint_{y_{\ell}=0}
\prod_{i=1}^{\ell}\frac{dy_i}{2\pi\i y_i}\Delta(y)\Big(\sum_{I_r}
T_{I_r}^{-1} A_{I_r}(y)
T_{I_r}\cdot\frac{T_{I_r}^{-1}\Delta(y)T_{I_r}}{\Delta(y)}\cdot
T^{-1}_{I_r}f(y^{-1})\Big) g(y). \ee Then the proof is provided by
the following Lemma $\Box$

\begin{lem} For any
$I_r=(i_1<i_2<\ldots<i_r)\subset\{1,2,\ldots,\ell\}$
the following relation holds.
$$
B_{I_r}(y)\,=\, \Big(\Delta_{I_r}(y)\cdot
T_{I_r}^{-1}A_{I_r}(y)T_{I_r}\Big)^{*},
$$
where
$$
\Delta_{I_r}(y)\,=\,\bigl(\Delta(y)\bigr)^{-1}\, T_{I_r}^{-1}\,
\Delta(y)\,T_{I_r},
$$
for all $i\in I_r$. Where given a  function $f(y)$ we define
$$
f\,^{*}(y)\,:=\,f(y^{-1})\,.
$$
\end{lem}
\noindent {\it Proof}: By  direct calculation one derives
$$
(T_{I_r}^{-1}A_{I_r}(y)T_{I_r})^{*}\,=\,\prod_{k=1}^r
\Big(1\,-\,q^{-1}\frac{y_{i_k-1}}{y_{i_k}}
\Big)^{1-\delta_{i_k-i_{k-1},\,1}},
$$
and
$$
\Delta_{I_r}^*(y)\,=\,\prod_{k=1}^r\,
\cfrac{\displaystyle\Big(1\,-\,\frac{y_{i_k}}{y_{i_k+1}}
\Big)^{1-\delta_{i_{k+1}-i_k,\,1}}}
{\displaystyle\Big(1\,-\,q^{-1}\frac{y_{i_k-1}}{y_{i_k}}
\Big)^{1-\delta_{i_k-i_{k-1},\,1}}},
$$
where we set
$i_0:=0,\,i_{r+1}:=\ell+1$. In this way we obtain
$$
(T_{I_r}^{-1}A_{I_r}(y)T_{I_r})^*\cdot\Delta_{I_r}^*(y)\,=\,
Y_{i_1}^{1-\delta_{i_2-i_1,\,1}}\cdot\ldots\cdot
Y_{i_{r-1}}^{1-\delta_{i_r-i_{r-1},\,1}}\cdot
Y_{i_r}^{1-\delta_{i_{r+1}-i_r,\,1}}$$
 $\Box$

To construct recursive formulas for $q$-deformed Whittaker functions
one should introduce a paring on a functions defined on the lattice
$\{y_i=q^{p_{\ell,i}+i-1}; i=1,\ldots,\ell\,;\, p_{\ell,i}\in
\IZ\,\}$ with appropriate decay at infinities. Let us define the
following analog of (\ref{todpair}) \be\label{latpair} \<f,g\>_{\rm
lat}\,=\,
\sum_{\underline{p}_{\ell}\in\ZZ^{\ell}}\Delta'(\underline{p}_{\ell})
f(-\underline{p}_{\ell})g(\underline{p}_{\ell}), \ee where
\be\label{latmeas}
\Delta'(\underline{p}_{\ell})=\prod_{i=1}^{\ell-1}
\Theta(p_{\ell,i+1}-p_{\ell,i})\,\, (p_{\ell,i+1}-p_{\ell,i})_q!
\,.\ee The following Proposition can be easily proved by mimicking
the proof of the Proposition \ref{selfadj}.
\begin{prop}\label{selfadlat}  Hamiltonians
$\CH_r^{\mathfrak{gl}_\ell}(\underline{p}_{\ell})$ and
$\widetilde{\CH}_r^{\mathfrak{gl}_\ell}(\underline{p}_{\ell})$ are adjoint with
respect to the pairing  (\ref{latpair})
\be\label{adlat}\<f,\CH_k^{\mathfrak{gl}_\ell} g\>_{\rm
lat}\,=\,\<\widetilde{\CH}_k^{\mathfrak{gl}_\ell} f, g\>_{\rm lat},\qquad
k=1,\ldots,\ell.\ee
\end{prop}

\subsection{Proof of Theorem \ref{mainTH}}

Now we are ready to prove Theorem \ref{mainTH}.  We use recursion
over the rank of $\mathfrak{gl}_k$. Set
$\,\Psi^{\mathfrak{gl}_{1}}_{\lambda_1}(p_{11})=q^{\l_1p_{11}}\,$
and assume that \be\label{eigprob}
\CH_r^{\mathfrak{gl}_{\ell}}(\underline{p}_{\ell})\cdot
\Psi^{\mathfrak{gl}_{\ell}}_{\l_1,\ldots ,\l_\ell}
(\underline{p}_{\ell})\,=\,
\chi^{\mathfrak{gl}_\ell}_r(q^{\sum_i\l_iE_{ii}})\,
\Psi^{\mathfrak{gl}_{\ell}}_{\l_1,\ldots,\l_{\ell}}
(\underline{p}_{\ell}),\\
\chi^{\mathfrak{gl}_\ell}_r(q^{\sum_i\l_iE_{ii}})\,=\, \sum_{I^{(\ell)}_r}
z_{i_1}z_{i_2}\cdot\ldots z_{i_r},\hspace{1.5cm}z_i=q^{\l_i}, \ee
where $I_r^{(\ell)}=\{i_1<i_2<\ldots<i_r\}\in (1,2,\ldots,\ell)$.

Let us define the function
$\Psi^{\mathfrak{gl}_{\ell+1}}_{\l_1,\ldots,\l_{\ell+1}}
(\underline{p}_{\ell+1})$ as follows \be\label{recprf}
\Psi^{\mathfrak{gl}_{\ell+1}}_{\l_1,\ldots,\l_{\ell+1}}
(\underline{p}_{\ell+1})=
\sum_{\underline{p}_{\ell}\in\ZZ^{\ell}}\,\,
\Delta'(\underline{p}_{\ell}) \,\, \CQ_{\ell+1,\ell}
(\underline{p}_{\ell+1},\underline{p}_{\ell})\\ \cdot
q^{\l_{\ell+1}(\sum_{i=1}^{\ell+1} p_{\ell+1,i}-\sum_{i=1}^{\ell}
p_{\ell,i})}\, \Psi^{\mathfrak{gl}_{\ell}}_{\l_1,\ldots,\l_{\ell}}
(\underline{p}_{\ell}),\ee where
$$
\CQ_{\ell+1,\ell}(\underline{p}_{\ell+1},\underline{p}_{\ell})\,=\,
\frac{ \prod_{i=1}^{\ell}\Theta(p_{\ell,i}-p_{\ell+1,i})
\Theta(p_{\ell+1,i+1}-p_{\ell,i})} {\prod_{i=1}^{\ell}
(p_{\ell,i}-p_{\ell+1,i})_q!\,\, (p_{\ell+1,i+1}-p_{\ell,i})_q!},
$$
and
$$ \Delta'(\underline{p}_{\ell})=
\prod_{i=1}^{\ell-1}\Theta(p_{\ell,i+1}-p_{\ell,i})\,\,
(p_{\ell,i+1}-p_{\ell,i})_q!\, .
$$
One should verify  the relations: \be
\CH_r^{\mathfrak{gl}_{\ell+1}}(\underline{p}_{\ell+1})\cdot
\Psi^{\mathfrak{gl}_{\ell+1}}_{\l_1,\ldots ,\l_{\ell+1}}
(\underline{p}_{\ell+1})\,=\,
\chi^{\mathfrak{gl}_{\ell+1}}_r(q^{\sum_i\l_iE_{ii}})\,
\Psi^{\mathfrak{gl}_{\ell+1}}_{\l_1,\ldots,\l_{\ell+1}}
(\underline{p}_{\ell+1}),\\
\chi^{\mathfrak{gl}_{\ell+1}}_r(q^{\sum_i\l_iE_{ii}})\,=\,
\sum_{I^{(\ell+1)}_r} z_{i_1}z_{i_2}\cdot\ldots
z_{i_r},\hspace{1.5cm}z_i=q^{\l_i}, \ee where
$I_r^{(\ell+1)}=\{i_1<i_2<\ldots<i_r\}\in (1,2,\ldots,\ell+1)$.

Applying Hamiltonians
$\CH_r^{\mathfrak{gl}_{\ell+1}}(\underline{p}_{\ell+1})$ to
(\ref{recprf}) and using intertwining relation given in Proposition
\ref{interlat}
  one can obtains
\be{{\CH}}_r^{\frak{gl}_{\ell+1}}(\underline{p}_{\ell+1})
\CQ_{\ell+1,\ell}
(\underline{p}_{\ell+1},\underline{p}_{\ell})\,\,\,
q^{\l_{\ell+1}(\sum_i p_{\ell+1,i}-\sum_i p_{\ell,i})}=\\=\nonumber
\Big\{q^{\l_{\ell+1}}
\widetilde{\CH}^{\frak{gl}_{\ell}}_{r-1}(-\underline{p}_{\ell})+
\widetilde{\CH}^{\frak{gl}_{\ell}}_r (-\underline{p}_{\ell})\Big\}
\CQ_{\ell+1,\ell}(\underline{p}_{\ell+1},\underline{p}_{\ell})\,\,\,
 q^{\l_{\ell+1}(\sum_i p_{\ell+1,i}-\sum_i p_{\ell,i})}.
\ee
 Now using  (\ref{adlat}), one obtains 
\be
\CH_r^{\mathfrak{gl}_{\ell+1}}(\underline{p}_{\ell+1})
\Psi^{\mathfrak{gl}_{\ell+1}}_{\l_1,\ldots ,\l_{\ell+1}}
(\underline{p}_{\ell+1})=\\ \nonumber= \sum_{\underline
{p}_{\ell}\in\ZZ^{\ell}}\,\, \Delta(\underline{p}_{\ell}) \,\,
\Big(\CH_r^{\mathfrak{gl}_{\ell}}(\underline{p}_{\ell})
\CQ_{\ell+1,\ell}(\underline{p}_{\ell+1},\underline{p}_{\ell}|q)\,
q^{\l_{\ell+1}(\sum_i p_{\ell+1,i}-\sum_i p_{\ell,i})}\Big)\,
\Psi^{\mathfrak{gl}_{\ell}}_{\l_1,\ldots,\l_{\ell}}(\underline{p}_{\ell})\\
\nonumber= \sum_{\underline{p}_{\ell}\in\ZZ^{\ell}}\,\,
\Delta(\underline{p}_{\ell}) \,\,(q^{\l_{\ell+1}}
\widetilde{\CH}^{\frak{gl}_{\ell}}_{r-1}(-\underline{p}_{\ell})+
\widetilde{\CH}^{\frak{gl}_{\ell}}_r
(-\underline{p}_{\ell}))\Big(\CQ_{\ell+1,\ell}
(\underline{p}_{\ell+1},\underline{p}_{\ell})\,
q^{\l_{\ell+1}(\sum_i p_{\ell+1,i}-\sum_i p_{\ell,i})}\Big)\,
\Psi^{\mathfrak{gl}_{\ell}}_{\l_1,\ldots,\l_{\ell}}
(\underline{p}_{\ell})\\
\nonumber=\sum_{\underline{p}_{\ell}\in\ZZ^{\ell}}\,\,
\Delta(\underline{p}_{\ell}) \,\, \Big(\CQ_{\ell+1,\ell}
(\underline{p}_{\ell+1},\underline{p}_{\ell})\,
q^{\l_{\ell+1}(\sum_i p_{\ell+1,i}-\sum_i p_{\ell,i})}\Big)\,
(q^{\l_{\ell+1}}{\CH}^{\frak{gl}_{\ell}}_{r-1}
(\underline{p}_{\ell})+{\CH}^{\frak{gl}_{\ell}}_r
(\underline{p}_{\ell}))
\Psi^{\mathfrak{gl}_{\ell}}_{\l_1,\ldots,\l_{\ell}}
(\underline{p}_{\ell})\\
\nonumber=\Big(q^{\l_{\ell+1}}\sum_{I^{(\ell)}_{r-1}}\prod_{i\in
I_{r-1}^{(\ell)}}q^{\l_i}+\sum_{I^{(\ell)}_{r}}\prod_{i\in
I_{r}^{(\ell)}}q^{\l_i}\Big)
\Psi^{\mathfrak{gl}_{\ell+1}}_{\l_1,\ldots,\l_{\ell+1}}
(\underline{p}_{\ell+1}) =(\sum_{ I^{(\ell+1)}_r}\prod\limits_{i\in
I^{(\ell+1)}_r}q^{\l_i})\,\,
\Psi^{\mathfrak{gl}_{\ell+1}}_{\l_1,\ldots ,\l_{\ell+1}}
(\underline{p}_{\ell+1}), \ee where
$I_r^{(\ell)}=\{i_1<i_2<\ldots<i_r\}\in (1,2,\ldots,\ell)$ and
$I_r^{(\ell+1)}=\{i_1<i_2<\ldots<i_r\}\in (1,2,\ldots,\ell+1)$. In
the last equality above we use the following relation.
$$
\chi_r^{\mathfrak{gl}_{\ell+1}}(\underline{z})\,=\,
z_{\ell+1}\chi^{\mathfrak{gl}_{\ell}}_{r-1}(\underline{z}')+
\chi_r^{\mathfrak{gl}_{\ell}}(\underline{z}'),
$$
where $\underline{z}'=(z_1,z_2,\ldots,z_\ell)$ for $z_i=q^{\l_i}$.
This completes the proof of Theorem \ref{mainTH}  $\Box$


\vskip 1cm

\noindent {\small {\bf A.G.}: {\sl Institute for Theoretical and
Experimental Physics, 117259, Moscow,  Russia; \hspace{8 cm}\,
\hphantom{xxx}  \hspace{4 mm} School of Mathematics, Trinity
College, Dublin 2, Ireland; \hspace{8 cm}\,
\hphantom{xxx}   \hspace{3 mm} Hamilton
Mathematics Institute, TCD, Dublin 2, Ireland;}}

\noindent{\small {\bf D.L.}: {\sl
 Institute for Theoretical and Experimental Physics,
117259, Moscow, Russia};\\
\hphantom{xxxxxx} {\it E-mail address}: {\tt lebedev@itep.ru}}\\

\noindent{\small {\bf S.O.}: {\sl
 Institute for Theoretical and Experimental Physics,
117259, Moscow, Russia};\\
\hphantom{xxxxxx} {\it E-mail address}: {\tt Sergey.Oblezin@itep.ru}}
\end{document}